\documentclass[11pt]{article}
\usepackage{amsfonts}

\usepackage{amsmath}

\newtheorem{theorem}{Theorem}

\newtheorem{lemma}[theorem]{Lemma}

\newtheorem{proposition}[theorem]{Proposition}

\newenvironment{proof}[1][Proof]{\textbf{#1.} }{\ \rule{0.5em}{0.5em}}

\begin{document}

\title{Integral representation of renormalized self-intersection local times}
\author{Yaozhong Hu\thanks{
Y. Hu is supported by the NSF grant DMS0504783.}, David Nualart\thanks{
D. Nualart is supported by the NSF grant DMS0604207.}, Jian Song \\
Department of Mathematics \\
University of Kansas \\
Lawrence, Kansas, 66045 USA}
\date{}
\maketitle

\begin{abstract}
In this paper we apply Clark-Ocone formula to deduce an explicit integral
representation for the renormalized self-intersection local time of the $d$%
-dimensional fractional Brownian motion with Hurst parameter $H\in (0,1)$.
  As a consequence, we derive the existence of some
exponential moments for this random variable.
\end{abstract}

\section{Introduction}

The purpose of this paper is to apply Clark-Ocone's formula to the
renormalized self-intersection local time of the $d$-dimensional fractional
Brownian motion. As a consequence, we derive the existence of some
exponential moments for this local time.

A well-known result in It\^{o}'s stochastic calculus asserts that any square
integrable random variable in the filtration generated by a $d$-dimensional
Brownian motion $W=\{W_{t},t\geq 0\}$ can be expressed as the sum of its
expectation plus \ the stochastic integral of a square integrable adapted
process:%
\begin{equation*}
F=E(F)+\sum_{i=1}^{d}\int_{0}^{\infty }u^{i}(t)dW_{t}^{i}.
\end{equation*}%
The process $u$ is determined by $F$, except on sets of measure zero. In
this context, Clark-Ocone formula provides an explicit representation of $u$
in terms of the derivative operator in the sense of Malliavin calculus. More
precisely, if $F$ belongs to the Sobolev space $\mathbb{D}^{1,2}$, then $%
u^{i}(t)=E(D_{t}^{i}F|\mathcal{F}_{t})$, where $D^{i}$ denotes the
derivative with respect to the $i$th component of the Brownian motion and $%
\left\{ \mathcal{F}_{t},t\geq 0\right\} $ is the filtration generated by the
Brownian motion. Extensions of this formula have been developed by \"{U}st%
\"{u}nel in \cite{Us}, and by Karatzas,\ Ocone and Li in \cite{KOL}.
Clark-Ocone formula has proved to be a useful tool in finding hedging
portfolios in mathematical finance (see, for instance, \cite{KO}).

The fractional Brownian motion on $\mathbb{R}^{d}$ with Hurst parameter \ $%
H\in (0,1)$ is a $d$-dimensional Gaussian process $B^{H}=\{B_{t}^{H},t\geq
0\}$ with zero mean and covariance function given by 
\begin{equation}
E(B_{t}^{H,i}B_{s}^{H,j})=\frac{\delta _{ij}}{2}(t^{2H}+s^{2H}-|t-s|^{2H}),
\label{x1}
\end{equation}%
where $i,j=1,\ldots ,d$,  $s,t\geq 0$, and 
\begin{equation*}
\delta _{ij}=\left\{ 
\begin{tabular}{lll}
$1$ & if & $i=j$ \\ 
$0$ &  & $i\neq j$%
\end{tabular}%
\right. 
\end{equation*}%
is the Kronecker symbol.  Assume $d\geq 2$. The \textit{self-intersection
local time }of $B^{H}$ is formally defined as%
\begin{equation*}
L=\int_{0}^{T}\int_{0}^{t}\delta _{0}(B_{t}^{H}-B_{s}^{H})ds,
\end{equation*}%
where $\delta _{0}$ is the Dirac delta function. It measures the amount of
time that the process spends intersecting itself on the time interval $[0,T]$%
. Rigorously, $L$ is defined as the limit in $L^{2}$, if it exists, of \ $%
L_{\varepsilon }=\int_{0}^{T}\int_{0}^{t}p_{\varepsilon
}(B_{t}^{H}-B_{s}^{H})dsdt$, as $\varepsilon $ tends to zero, where $%
p_{\varepsilon }$ denotes the heat kernel.

For $H=\frac{1}{2}$, the process $B^{H}$ is a classical Brownian motion and
its self-intersection local time has been studied by many authors (see
Albeverio et al. \cite{AHZ}, Calais and Yor \cite{CY}, He et al. \cite{HYYW}%
, Hu \cite{Hu}, Imkeller et al. \cite{IPV}, Varadhan \cite{Var}, Yor\cite%
{Yor}, and the references therein). \ In this case, if $d=2$, Varadhan \cite%
{Var} has proved that \ $L_{\varepsilon }$ does not converge in $L^{2}$, but
it can be renormalized so that $L_{\varepsilon }-E(L_{\varepsilon })$
converges in $L^{2}$ as $\varepsilon $ tends to zero to a random variable
that we denote by $\widetilde{L}$. This result has been extended by Rosen %
\cite{Ro} to the case $\ $ $H\in \left( \frac{1}{2},\frac{3}{4}\right) $
(still when $d=2$), and by Hu and Nualart in \cite{HN}, where they have
obtained the following complete result on the existence of the
self-intersection local time of the fractional Brownian motion:

\begin{itemize}
\item[(i)] The self-intersection local time $L$ exists if and only if $Hd<1$.

\item[(ii)] If $Hd\geq 1$, the renormalized self-intersection local time $%
\widetilde{L}$ exists if and only if $Hd<\frac{3}{2}$.
\end{itemize}

An important question is the existence of moments and exponential moments
for the (renormalized) self-intersection local time. Along this direction,\
Le Gall \cite{LG} proved that for the planar Brownian motion, there is a
critical exponent $\lambda _{0}$, such that $E\left( \exp \lambda \widetilde{%
L}\right) <\infty $ for all $\lambda <\lambda _{0}$, and $E\left( \exp
\lambda \widetilde{L}\right) =\infty $ if $\lambda >\lambda _{0}$. Using the
theory of large deviations, Bass and Chen proved in \ \cite{BC} that the
critical exponent $\lambda _{0}$ coincides with $A^{-4}$, where $A$ is the \
best constant in the Gagliardo-Nirenberg inequality.

Clark-Ocone formula seems to be a suitable tool to analyze the renormalized
self-intersection local time, because in this formula we do not take into
account the expectation of the random variable. The fractional Brownian
motion can be expressed as the stochastic integral 
\begin{equation*}
B^H_t= \int_0^t K_H(t,s) dW_s
\end{equation*}
of a square integrable kernel $K_H(t,s)$ with respect to an underlying
Brownian motion $W$. In this way the renormalized self-intersection local
time $\widetilde{L}$ is a functional of the Brownian motion $W$, and we can
obtain an explicit integral representation $\widetilde{L}$, in the general
case $Hd<\frac{3}{2}$. This formula allows us to obtain some exponential
moments for the renormalized self-intersection local time, using the method
of moments.

The paper is organized as follows. In Section 2 we present some
preliminaries on Malliavin calculus and Clark-Ocone formula. Section 3 is
devoted to derive estimates for the moments of the self-intersection local
time in the case of a general $d$-dimensional Gaussian process, using the
method of moments. In the case of the fractional Brownian motion, this
provides the existence of exponential moments in the case $Hd<1$. Section 4
contains the main result, which is the integral representation of the
renormalized self-intersection local time of the fractional Brownian motion
in the case $H<\min \left( \frac{3}{2d},\frac{2}{d+1}\right) $. As an
application we show that $E\left( \exp \left| \widetilde{L}\right|
^{p}\right) <\infty $ if $p<\frac{1}{2}\left[ \left( \frac{1}{2}+H\right)
\left( \frac{d}{2}-\frac{1}{4H}\right) \right] ^{-1}$. A crucial tool is the
local nondeterminism property introduced by Berman \ in {\cite{Be}} and
developed by many authors (see Xiao \cite{Xiao} and the references therein).

\setcounter{equation}{0}

\section{Preliminaries on Malliavin calculus and \protect\break Clark-Ocone
formula}

We need some preliminaries on the Malliavin calculus for the $d$-dimensional
Brownian motion $W=\{W_{t},t\geq 0\}$. We refer to Malliavin \cite{Ma} and
Nualart \cite{Nu} for a more detailed presentation of this theory.

We assume that $W$ is defined in a complete probability space $(\Omega ,%
\mathcal{F},P)$, and the $\sigma $-field $\mathcal{F}$ is generated by $W$.
Let us denote by $H$ the Hilbert space $L^{2}(\mathbb{R}_{+};\mathbb{R}^{d})$%
, and for any function $h\in H$ we set%
\begin{equation*}
W(h)=\sum_{i=1}^{d}\int_{0}^{\infty }h^{i}(t)dW_{t}^{i}.
\end{equation*}%
Let $\mathcal{S}$ be the class of smooth and cylindrical random variables of
the form%
\begin{equation*}
F=f(W(h_{1}),\ldots ,W(h_{n})),
\end{equation*}%
where $n\geq 1$, $h_{1},\ldots ,h_{n}\in H$, and $f$ is an infinitely
differentiable function such that together with all its partial derivatives
has at most polynomial growth order. The derivative operator of the random
variable $F$ is defined as%
\begin{equation*}
D_{t}^{i}F=\sum_{j=1}^{n}\frac{\partial f}{\partial x_{j}}(W(h_{1}),\ldots
,W(h_{n}))h_{j}^{i}(t),
\end{equation*}%
where $i=1,\ldots ,d$ and $t\geq 0$. In this way, we interpret $DF$ as a
random variable with values in the Hilbert space $H$. The derivative is a
closable operator on $L^{2}(\Omega )$ with values in $L^{2}(\Omega ;H)$. We
denote by $\mathbb{D}^{1,2}$ the Hilbert spaced defined as the completion of 
$\mathcal{S}$ with respect to the scalar product%
\begin{equation*}
\left\langle F,G\right\rangle _{1,2}=E(FG)+E\left(
\sum_{i=1}^{d}\int_{0}^{\infty }D_{t}^{i}FD_{t}^{i}Gdt\right) .
\end{equation*}%
The divergence operator $\delta $ is the adjoint of the derivative operator $%
D$. The operator $\delta $ is an unbounded operator from \ $L^{2}(\Omega ;H)$
into $L^{2}(\Omega )$, and is determined by the duality relationship%
\begin{equation*}
E(\delta (u)F)=E(\left\langle u,DF\right\rangle _{H}),
\end{equation*}%
for any $u$ in the domain of $\delta $, and $F$ in \ $\mathbb{D}^{1,2}$.
Gaveau and Trauber \cite{GT} proved that $\delta $ is an extension of the
classical It\^{o} integral in the sense that any $d$-dimensional square
integrable adapted process belongs to the domain of $\delta $, and $\delta
(u)$ coincides with the It\^{o} integral of $u$:%
\begin{equation*}
\delta (u)=\sum_{i=1}^{d}\int_{0}^{\infty }u^{i}(t)dW_{t}^{i}.
\end{equation*}%
It is well-known that any random variable $F\in L^{2}(\Omega )$, possesses a
stochastic integral representation of the form%
\begin{equation*}
F=E(F)+\sum_{i=1}^{d}\int_{0}^{\infty }u^{i}(t)dW_{t}^{i},
\end{equation*}%
for some $d$-dimensional square integrable adapted process $u$. Clark-Ocone
formula says that if $F\in \mathbb{D}^{1,2}$, then%
\begin{equation}
F=E(F)+\sum_{i=1}^{d}\int_{0}^{\infty }E(D_{t}^{i}F|\mathcal{F}%
_{t})dW_{t}^{i}.  \label{CO}
\end{equation}

\setcounter{equation}{0}

\section{Exponential integrability of the self-intersection local time}

Suppose that $W=\{W_{t},t\geq 0\}$ is a \ $d$-dimensional standard $\ $
Brownian motion, defined in a complete probability space $(\Omega ,\mathcal{F%
},P)$. Suppose that $\mathcal{F}$ is generated by $W$. We denote by $\left\{ 
\mathcal{F}_{t},t\geq 0\right\} $ the filtration generated by $W$ and the
sets of probability zero. Consider a $d$-dimensional Gaussian processs of
the form%
\begin{equation}
B_{t}=\int_{0}^{t}K(t,s)dW_{s},  \label{g1}
\end{equation}%
where $K(t,s)$ is a measurable kernel satisfying \ $%
\int_{0}^{t}K(t,s)^{2}ds<\infty $ for all $t\geq 0$. We will assume that $%
K(t,s)=0$ if $s>t$. \ 

Fix a time interval $[0,T]$. We will make use of the following property on
the kernel $K(t,s)$:

\textbf{(H1) }\textit{For any }$s,t\in \lbrack 0,T]$\textit{, }$s<t$\textit{%
\ we have}%
\begin{equation}
\int_{s}^{t}K(t,\theta )^{2}d\theta \geq k_{1}(t-s)^{2H}  \label{c1}
\end{equation}%
\textit{for some constants }$k_{1}>0$\textit{, and} $H\in (0,1)$.

Notice that $\mathrm{Var}\left( B_{t}^{i}|\mathcal{F}_{s}\right)
=\int_{s}^{t}K(t,\theta )^{2}d\theta $, so condition \textbf{(H1)} is
equivalent to say that $\mathrm{Var}\left( B_{t}^{i}|\mathcal{F}_{s}\right)
\geq k_{1}(t-s)^{2H}$, for each component $i=1,\ldots ,d$. This property is
satisfied, for instance, in the following two examples:

\medskip

\noindent \textbf{Example 1} Suppose that $K(t,s)=(t-s)^{H-\frac{1}{2}}$ .
Then, we have equality in \ (\ref{c1}) with $k_{1}=\frac{1}{2H}$.

\medskip

\noindent \textbf{Example 2 }Condition \textbf{(H1)} is satisfied by the
kernel of the fractional Brownian motion, as a consequence of the local
nondeterminism property (see \ (\ref{LND}) below).

\medskip

We will denote by $C$ a generic constant depending on $T$, the dimension\ $d$%
, and the constants appearing in the hypothesis such as $H$ and $k_{1}$.

The \textit{self-intersection local time} of the process $B$ in the time
interval $[0,T]$, denoted by $L$, is defined as the limit in $L^{2}$ as $%
\varepsilon $ tends to zero of%
\begin{equation}
L_{\varepsilon }=\int_{0}^{T}\int_{0}^{t}p_{\varepsilon }(B_{t}-B_{s})ds,
\label{e1}
\end{equation}%
where $p_{\varepsilon }$ denotes the heat kernel 
\begin{equation*}
p_{\varepsilon }(x)=(2\pi \varepsilon )^{-\frac{d}{2}}\exp \left( -\frac{%
|x|^{2}}{2\varepsilon }\right) .
\end{equation*}%
The next theorem asserts that $L$ exists if $Hd<1$, and it has exponential
moments of order $\frac{1}{Hd}$.

\begin{theorem}
Suppose that $Hd<1$. Then, the self-intersection local time $L$ exists as
the limit in $L^{2}$ of $L_{\varepsilon }$, as $\varepsilon $ tends to zero,
and for all integers $n\geq 1$ we have%
\begin{equation*}
E(L^{n})\leq C^{n}\left( n!\right) ^{Hd},
\end{equation*}%
for some constant $C$. As a consequence, $\ $%
\begin{equation*}
E(e^{L^{p}})<\infty ,
\end{equation*}%
for any $p<\frac{1}{Hd}$, and there exists a constant $\lambda_0>0$ such that
$E(e^{\lambda L^{\frac 1{Hd}}})<\infty$ for all $\lambda<\lambda_0$.
\end{theorem}

\begin{proof}
From the equality%
\begin{equation*}
p_{\varepsilon }(x)=\frac{1}{(2\pi )^{d}}\int_{\mathbb{R}^{d}}\exp \left(
i\left\langle \xi ,x\right\rangle -\frac{\varepsilon |\xi |^{2}}{2}\right)
d\xi
\end{equation*}%
and the definition of $L_{\varepsilon }$, we obtain%
\begin{equation*}
L_{\varepsilon }=\frac{1}{(2\pi )^{d}}\int_{0}^{T}\int_{0}^{t}\int_{\mathbb{R%
}^{d}}\exp \left( i\left\langle \xi ,B_{t}-B_{s}\right\rangle -\frac{%
\varepsilon |\xi |^{2}}{2}\right) d\xi dsdt.
\end{equation*}%
This expression allows us to compute the moments of $L_{\varepsilon }$. Fix
an integer $n\geq 1$. Denote by $T_{n}$ the set $\{0<s<t<T\}^{n}$. Then%
\begin{eqnarray}
E(L_{\varepsilon }^{n}) &=&\frac{1}{(2\pi )^{nd}}\int_{T_{n}}\int_{\mathbb{R}%
^{nd}}E\left[ \exp \left( i\left\langle \xi
_{1},B_{t_{1}}-B_{s_{1}}\right\rangle +\cdots +i\left\langle \xi
_{n},B_{t_{n}}-B_{s_{n}}\right\rangle \right) \right]  \notag \\
&&\times \exp \left( -\frac{\varepsilon }{2}\sum_{j=1}^{n}|\xi
_{j}|^{2}\right) d\xi _{1}\cdots d\xi _{n}dsdt\mathbf{,}  \label{f2}
\end{eqnarray}%
where $s=(s_{1},\ldots ,s_{n})$ and \ $t=(t_{1},\ldots ,t_{n})$. Notice that%
\begin{eqnarray}
&&\int_{\mathbb{R}^{nd}}E\left[ \exp \left( i\left\langle \xi
_{1},B_{t_{1}}-B_{s_{1}}\right\rangle +\cdots +i\left\langle \xi
_{n},B_{t_{n}}-B_{s_{n}}\right\rangle \right) \right]  \notag \\
&&\times e^{-\frac{\varepsilon }{2}\sum_{j=1}^{n}|\xi _{j}|^{2}}d\xi
_{1}\cdots d\xi _{n}  \notag \\
&=&\int_{\mathbb{R}^{nd}}\exp \left( -\frac{1}{2}E\left[ \left( \left\langle
\xi _{1},B_{t_{1}}-B_{s_{1}}\right\rangle +\cdots +\left\langle \xi
_{n},B_{t_{n}}-B_{s_{n}}\right\rangle \right) ^{2}\right] \right)  \notag \\
&&\times e^{-\frac{\varepsilon }{2}\sum_{j=1}^{n}|\xi _{j}|^{2}}d\xi
_{1}\cdots d\xi _{n}  \notag \\
&=&\left( \int_{\mathbb{R}^{n}}\exp \left( -\frac{1}{2}\xi ^{T}Q\xi \right)
e^{-\frac{\varepsilon }{2}|\xi |^{2}}d\xi \right) ^{d},  \label{f1}
\end{eqnarray}%
where $Q$ is the covariance matrix of the $n$-dimensional random vector $%
(B_{t_{1}}^{1}-B_{s_{1}}^{1},\ldots ,B_{t_{n}}^{1}-B_{s_{n}}^{1})$.
Substituting (\ref{f1}) into (\ref{f2}) yields%
\begin{equation*}
E(L_{\varepsilon }^{n})=\frac{1}{(2\pi )^{nd}}\int_{T_{n}}\left( \int_{%
\mathbb{R}^{n}}\exp \left( -\frac{1}{2}\xi ^{T}Q\xi \right) e^{-\frac{%
\varepsilon }{2}|\xi |^{2}}d\xi \right) ^{d}dsdt,
\end{equation*}%
and $E(L_{\varepsilon }^{n})$ converges as $\varepsilon $ tends to zero to%
\begin{eqnarray*}
\alpha _{n} &=&\frac{1}{(2\pi )^{nd}}\int_{T_{n}}\left( \int_{\mathbb{R}%
^{n}}\exp \left( -\frac{1}{2}\xi ^{T}Q\xi \right) d\xi \right) ^{d}dsdt \\
&=&\frac{1}{(2\pi )^{\frac{nd}{2}}}\int_{T_{n}}\left( \det Q\right) ^{-\frac{%
d}{2}}dsdt,
\end{eqnarray*}%
provided $\alpha _{n}$ is finite.

If $\alpha_2<\infty$, then in the same way as before we obtain 
\begin{equation*}
\lim_{\varepsilon , \delta \downarrow 0} E(L_\varepsilon L_\delta)=\alpha_2,
\end{equation*}
which implies that $L_\varepsilon$ converges in $L^2$ as $\varepsilon$ tends
to zero. Furthermore, if $\alpha _{n}$ is finite for all $n\geq 1$, then we
deduce the convergence in $L^{p}$ for any $p\geq 2$ of $L_{\varepsilon }$ as 
$\varepsilon $ tends to zero. The limit, denoted by $L$, will be, by
definition, the self-intersection local time of the process $B$ in the time
interval $[0,T]$. To complete the proof of the theorem it suffices to show
that $\alpha _{n}$ is bounded by $C^{n}\left( n!\right) ^{Hd}$, for some
constant $C$.

We can write%
\begin{equation*}
\alpha _{n}=\frac{n!}{\ (2\pi )^{\frac{nd}{2}}}\int_{T_{n}\cap
\{t_{1}<\cdots <t_{n}\}}\left( \det Q\right) ^{-\frac{d}{2}}dsdt\mathbf{.}
\end{equation*}%
For each $i=1,\ldots ,n$ we denote by $\tau _{i}$ the point in the set $%
\{s_{i},s_{i+1},\ldots ,s_{n},t_{i-1}\}$ which is closer to $t_{i}$ from the
left. Then, by \textbf{(H1)} and the fact that $s_{i}<t_{i}$, $i=1,\dots ,n$%
, we obtain, using Lemma \ref{lem2} in the Appendix, 
\begin{eqnarray*}
\det Q &=&\mathrm{Var}(B_{t_{1}}^{1}-B_{s_{1}}^{1})\mathrm{Var}%
(B_{t_{2}}^{1}-B_{s_{2}}^{1}|B_{t_{1}}^{1}-B_{s_{1}}^{1}) \\
&&\times \cdots \times \mathrm{Var}%
(B_{t_{n}}^{1}-B_{s_{n}}^{1}|B_{t_{1}}^{1}-B_{s_{1}}^{1},\dots
,B_{t_{n-1}}^{1}-B_{s_{n-1}}^{1}) \\
&\geq &\mathrm{Var}(B_{t_{1}}^{1}|B_{s_{1}}^{1})\mathrm{Var}%
(B_{t_{2}}^{1}|B_{t_{1}}^{1},B_{s_{1}}^{1},B_{s_{2}}^{1}) \\
&&\times \cdots \times \mathrm{Var}%
(B_{t_{n}}^{1}|B_{t_{1}}^{1},B_{s_{1}}^{1},\dots
,B_{t_{n-1}}^{1},B_{s_{n-1}}^{1},B_{s_{n}}^{1}) \\
&\geq &\mathrm{Var}(B_{t_{1}}^{1}|\mathcal{F}_{\tau _{1}})\mathrm{Var}%
(B_{t_{2}}^{1}|\mathcal{F}_{\tau _{2}})\cdots \mathrm{Var}(B_{t_{n}}^{1}|%
\mathcal{F}_{\tau _{n}}) \\
&\geq &k_{1}^{n}(t_{1}-\tau _{1})^{2H}(t_{2}-\tau _{2})^{2H}\cdots
(t_{n}-\tau _{n})^{2H}.
\end{eqnarray*}%
As a consequence,%
\begin{equation*}
\alpha _{n}\leq \frac{n!}{\ (2\pi )^{\frac{nd}{2}}}k_{1}^{-\frac{nd}{2}%
}\int_{T_{n}\cap \{t_{1}<\cdots <t_{n}\}}\prod_{i=1}^{n}(t_{i}-\tau
_{i})^{-Hd}dsdt.
\end{equation*}%
If we fix the points $t_{1}<\cdots <t_{n}$, there are $3\times 5\times
\cdots \times (2n-1)=(2n-1)!!$ posible ways to place the points $%
s_{1},\ldots ,s_{n}$. In fact, $s_{1}$ must be in $(0,t_{1})$. For $s_{2}$
we have three choices: $(0,s_{1})$, $(s_{1},t_{1})$ and $(t_{1},t_{2})$. By
a recursive argument it is clear that we have $(2i-1)$ possible choices for $%
s_{i}$, given $s_{1},\ldots ,s_{i-1}$. In this way, up to a set of measure
zero, we can decompose the set $T_{n}\cap \{t_{1}<\cdots <t_{n}\}$ into the
union of $(2n-1)!!$ disjoint subsets. The integral of $%
\prod_{i=1}^{n}(t_{i}-\tau _{i})^{-Hd}$ on each one of these subset can be
expressed as%
\begin{equation*}
\Phi _{\sigma }=\int_{\{0<z_{1}<\cdots <z_{2n}<T\}}\prod_{i=1}^{n}(z_{\sigma
(i)}-z_{\sigma (i)-1})^{-Hd}dz,
\end{equation*}%
where $\sigma (1)<\cdots <\sigma (n)$ are $n$ elements \ in $\{1,2,\ldots
,2n\}$, and $z=(z_{1},\dots ,z_{2n})$. Making the change of variables $%
y_{i}=z_{i}-z_{i-1}$, $i=1,\ldots ,2n$ (with the convention $z_{0}=0$) we
obtain%
\begin{eqnarray*}
\Phi _{\sigma } &=&\int_{\{0<y_{1}+\cdots
+y_{2n}<T\}}\prod_{i=1}^{n}y_{\sigma (i)}^{-Hd}dy\leq \frac{T^{n}}{n!}%
\int_{\{0<y_{1}+\cdots +y_{n}<T\}}\prod_{i=1}^{n}y_{i}^{-Hd}dy \\
&=&\frac{1}{n!}T^{n(2-Hd)+Hd}\frac{\Gamma (1-Hd)^{n-1}}{\Gamma (n(1-Hd)+Hd+1)%
}.
\end{eqnarray*}%
Therefore%
\begin{eqnarray*}
\alpha _{n} &\leq &\frac{k_{1}^{-\frac{nd}{2}}(2n-1)!!T^{n(2-Hd)+Hd}\Gamma
(1-Hd)^{n-1}}{\ (2\pi )^{\frac{nd}{2}}\Gamma (n(1-Hd)+Hd+1)} \\
&=&C_{1}C_{2}^{n}\frac{(2n-1)!!}{\Gamma (n(1-Hd)+Hd+1)},
\end{eqnarray*}%
with $C_{1}=T^{Hd}\Gamma (1-Hd)^{-1}$, and $C_{2}=\frac{k_{1}^{-\frac{d}{2}%
}\Gamma (1-Hd)T^{2-Hd}}{(2\pi )^{\frac{d}{2}}}$. Taking into account that $%
(2n-1)!!\leq 2^{n-1}n!$, and that%
\begin{equation*}
\Gamma (n(1-Hd)+Hd+1)\geq C^{n}(n!)^{1-Hd},
\end{equation*}%
for some constant $C$, we obtain the desired estimate.
\end{proof}

If $Hd\geq 1$, the above result is no longer true. In that case the
expectation of $L_{\varepsilon }$ blows up as $\varepsilon $ tends to zero.
In fact, if we denote   $\sigma ^{2}(s,t)=\mathrm{Var}(B_{t}^{1}-B_{s}^{1})$%
, for $s<t$, then 
\begin{equation*}
E(L_{\varepsilon })=\int_{0}^{T}\int_{0}^{t}p_{\varepsilon +\sigma
^{2}(s,t)}(0)dsdt=(2\pi )^{-\frac{d}{2}}\int_{0}^{T}\int_{0}^{t}(\varepsilon
+\sigma ^{2}(s,t))^{-\frac{d}{2}}dsdt,
\end{equation*}%
which converges to%
\begin{equation*}
(2\pi )^{-\frac{d}{2}}\int_{0}^{T}\int_{0}^{t}\sigma ^{2}(s,t)^{-\frac{d}{2}%
}dsdt\geq (2\pi )^{-\frac{d}{2}}k_{1}^{-\frac{d}{2}}\int_{0}^{T}%
\int_{0}^{t}(t-s)^{-Hd}dsdt=\infty .
\end{equation*}%
In this case, one can study the existence of the renormalized
self-intersection local time defined as the limit as $\varepsilon $ tends to
zero of $L_{\varepsilon }-E(L_{\varepsilon })$. In the next section we
discuss the existence and exponential moments of the renormalized
self-intersection local time, using Clark-Ocone formula, in the case of the
fractional Brownian motion.

\setcounter{equation}{0}

\section{Renormalized self-intersection local time of the fBm}

The fractional Brownian motion on $\mathbb{R}^{d}$ with Hurst parameter $%
H\in (0,1)$ is a $d$-dimensional Gaussian process $B^{H}=\{B_{t}^{H},t\geq
0\}$ with zero mean and covariance function given by (\ref{x1}). We will
assume that $d\geq 2$.

It is well-known that $B^{H}$ possesses the following integral representation%
\begin{equation*}
B_{t}^{H}=\int_{0}^{t}K_{H}(t,s)dW_{s},
\end{equation*}%
where $W=\{W_{t},t\geq 0\}$ is a $d$-dimensional Brownian motion, and $%
K_{H}(s,t)$ is the square integrable kernel given by%
\begin{equation*}
K_{H}(t,s)=C_{H,1}s^{\frac{1}{2}-H}\int_{s}^{t}(u-s)^{H-\frac{3}{2}}u^{H-%
\frac{1}{2}}du,\newline
\end{equation*}%
if $H>\frac{1}{2}$, and by%
\begin{equation*}
K_{H}(t,s)=C_{H,2}\left[ \left( \frac{t}{s}\right) ^{H-\frac{1}{2}}(t-s)^{H-%
\frac{1}{2}}-(H-\frac{1}{2})s^{\frac{1}{2}-H}\int_{s}^{t}u^{H-\frac{3}{2}%
}(u-s)^{H-\frac{1}{2}}du\right] ,\newline
\end{equation*}%
if $H<\frac{1}{2}$, for any $s<t$, where the constants are $C_{H,1}=\left[ 
\frac{H(2H-1)}{B(2-2H,H-\frac{1}{2})}\right] ^{\frac{1}{2}}$, and $C_{H,2}=%
\left[ \frac{2H}{(1-2H)b(1-2H,H+\frac{1}{2})}\right] ^{\frac{1}{2}}$, where $%
B(\alpha ,\beta )$ denotes th beta function.

The processes $B^{H}$ and $W$ generate the same filtration, that is, $%
\mathcal{F}_{t}=\sigma \{W_{s},0\leq s\leq t\}=\sigma \{B_{s}^{H},0\leq
s\leq t\}$.

The \ fractional Brownian motion satisfies the following local
nondeterminism property:

\medskip

\textbf{(LND)} \textit{There exists a constant }$k_{2}>0$\textit{, depending
only on }$H$\textit{\ and }$T$\textit{, such that for any }$t\in \lbrack
0,T] $, $0<r<t\wedge (T-t)$ and for $i=1,\ldots ,d$,

\begin{equation}
\mathrm{Var}(B_{t}^{H,i}|B_{s}^{H,i}:|s-t|\geq r)\ \geq k_{2\ }r^{2H}.
\label{LND}
\end{equation}

Consider the approximated self-intersection local time $L_{\varepsilon }$
introduced in (\ref{e1}). From the general result proved in Section 2 it
follows that if $Hd<1$, then $L_{\varepsilon }$ converges in $L^{2}$ to the
self-intersection local time $L$, and the random variable $L$ has
exponential moments. If $Hd\geq 1$, this result is no longer true, and one
considers the renormalization of the self-intersection local time,
introduced by Varadhan.

The purpose of this section is to apply the Clark-Ocone formula to provide a
stochastic integral representation for the renormalized self-intersection
local time $\widetilde{L}$. As a consequence, we will prove the existence of
some exponential moments for the random variable \ $\widetilde{L}$.

\bigskip

\begin{theorem}
\label{thm1}Suppose that $H<\min \left( \frac{3}{2d},\frac{2}{d+1}\right) $.
Then the renormalized self-intersection local time of the $d$-dimensional
fractional Brownian motion $B^{H}$ exists in $L^{2}$ and it has the
following integral representation%
\begin{equation}
\widetilde{L}=-\sum_{i=1}^{d}\int_{0}^{T}\left( \int_{r}^{T}\int_{0}^{t}%
\frac{A_{r,t,s}^{i}}{\sigma _{r,s,t}^{2}}p_{\sigma
_{r,s,t}^{2}}(A_{r,t,s}^{i})\left[ K_{H}(t,r)-K_{H}(s,r)\right] dsdt\right)
dW_{r}^{i},  \label{h1}
\end{equation}%
where%
\begin{equation*}
A_{r,t,s}=E(B_{t}^{H}-B_{s}^{H}|\mathcal{F}_{r})
\end{equation*}%
and%
\begin{equation*}
\sigma _{r,s,t}^{2}=\mathrm{Var}(B_{t}^{H,i}-B_{s}^{H,i}|\mathcal{F}_{r}).
\end{equation*}
\end{theorem}

\begin{proof}
The proof will be done in several steps.

\textbf{Step 1} We are going to apply Clark-Ocone formula to the random
variable $L_{\varepsilon }$. It is clear that $L_{\varepsilon }$ belongs to $%
\mathbb{D}^{1,2}$, and its derivative can be computed as follows%
\begin{equation*}
D_{r}^{i}L_{\varepsilon }=\int_{0}^{T}\int_{0}^{t}\frac{\partial
p_{\varepsilon }}{\partial x_{i}}(B_{t}^{H}-B_{s}^{H})D_{r}^{i}\left(
B_{t}^{H,i}-B_{s}^{H,i}\right) dsdt,
\end{equation*}%
where $r\in \lbrack 0,T]$, and $i=1,\ldots ,d$. Using 
\begin{equation*}
D_{r}^{i}\left( B_{t}^{H,i}-B_{s}^{H,i}\right) =\left[ K_{H}(t,r)-K_{H}(s,r)%
\right] \mathbf{1}_{[0,t]}(r),
\end{equation*}%
we obtain%
\begin{equation}
D_{r}^{i}L_{\varepsilon }=\int_{r}^{T}\int_{0}^{t}\frac{\partial
p_{\varepsilon }}{\partial x_{i}}(B_{t}^{H}-B_{s}^{H})\left[
K_{H}(t,r)-K_{H}(s,r)\right] dsdt.  \label{e2}
\end{equation}%
The next step is to compute the conditional expectation $E(D_{r}^{i}L_{%
\varepsilon }|\mathcal{F}_{r})$. The conditional law of $B_{t}^{H}-B_{s}^{H}$
given $\mathcal{F}_{r}$ is normal with mean $A_{r,t,s}$ and covariance
matrix $\sigma _{r,s,t}^{2}I_{d}$, where $I_{d}$ is the $d$-dimensional
identity matrix. Hence, the conditional expectation $E\left( \frac{\partial
p_{\varepsilon }}{\partial x_{i}}(B_{t}^{H}-B_{s}^{H})|\mathcal{F}%
_{r}\right) $ is given by%
\begin{eqnarray*}
E\left( \frac{\partial p_{\varepsilon }}{\partial x_{i}}%
(B_{t}^{H}-B_{s}^{H})|\mathcal{F}_{r}\right)  &=&\int_{\mathbb{R}^{d}}\frac{%
\partial p_{\varepsilon }}{\partial x_{i}}(y)p_{\sigma
_{r,s,t}^{2}}(y-A_{r,t,s})dy \\
&=&\frac{\partial p_{\varepsilon +\sigma _{r,s,t}^{2}}}{\partial x_{i}}%
(A_{r,t,s}) \\
&=&-\frac{A_{r,t,s}^{i}}{\varepsilon +\sigma _{r,s,t}^{2}}p_{\varepsilon
+\sigma _{r,s,t}^{2}}(A_{r,t,s}).
\end{eqnarray*}%
As a consequence, from (\ref{e2}) we obtain%
\begin{equation*}
E\left( D_{r}^{i}L_{\varepsilon }|\mathcal{F}_{r}\right)
=-\int_{r}^{T}\int_{0}^{t}\frac{A_{r,t,s}^{i}}{\varepsilon +\sigma
_{r,s,t}^{2}}p_{\varepsilon +\sigma _{r,s,t}^{2}}(A_{r,t,s})\left[
K_{H}(t,r)-K_{H}(s,r)\right] dsdt,
\end{equation*}%
and this leads to the following integral representation for $L_{\varepsilon
}-E(L_{\varepsilon })$%
\begin{eqnarray*}
&&L_{\varepsilon }-E(L_{\varepsilon }) \\
&=&-\sum_{i=1}^{d}\int_{0}^{T}\left( \int_{r}^{T}\int_{0}^{t}\frac{%
A_{r,t,s}^{i}}{\varepsilon +\sigma _{r,s,t}^{2}}p_{\varepsilon +\sigma
_{r,s,t}^{2}}(A_{r,t,s})\left[ K_{H}(t,r)-K_{H}(s,r)\right] dsdt\right)
dW_{r}^{i}.
\end{eqnarray*}

\textbf{Step 2 }In order to pass to the limit as $\varepsilon $ tends to
zero we proceed as follows. Set%
\begin{equation}
\Sigma _{\varepsilon }^{i}(r,t,s)=\frac{A_{r,t,s}^{i}}{\varepsilon +\sigma
_{r,s,t}^{2}}p_{\varepsilon +\sigma _{r,s,t}^{2}}(A_{r,t,s})\left[
K_{H}(t,r)-K_{H}(s,r)\right] .  \label{h5}
\end{equation}%
Clearly, $\Sigma _{\varepsilon }^{i}(r,t,s)$ converges pointwise as $%
\varepsilon $ tends to zero to 
\begin{equation*}
\Sigma ^{i}(r,t,s)=\frac{A_{r,t,s}^{i}}{\sigma _{r,s,t}^{2}}p_{\sigma
_{r,s,t}^{2}}(A_{r,t,s})\left[ K_{H}(t,r)-K_{H}(s,r)\right] .
\end{equation*}%
In order to establish the convergence of the integrals in the variables $s$
and $t$, we will first decompose the interval $\ \ [0,t]\ $ into the
disjoint union of $[r,t]$ and $[0,r)$. In this way we obtain%
\begin{equation*}
L_{\varepsilon }-E(L_{\varepsilon })=L_{\varepsilon }^{(1)}+L_{\varepsilon
}^{(2)},
\end{equation*}%
where%
\begin{equation*}
L_{\varepsilon }^{(1)}=-\sum_{i=1}^{d}\int_{0}^{T}\left(
\int_{r}^{T}\int_{r}^{t}\Sigma _{\varepsilon }^{i}(r,t,s)dsdt\right)
dW_{r}^{i},
\end{equation*}%
and%
\begin{equation*}
L_{\varepsilon }^{(2)}=-\sum_{i=1}^{d}\int_{0}^{T}\left(
\int_{r}^{T}\int_{0}^{r}\Sigma _{\varepsilon }^{i}(r,t,s)dsdt\right)
dW_{r}^{i}.
\end{equation*}

\textbf{Step 3 }We claim that the random field $\Sigma _{\varepsilon
}^{i}(r,t,s)$ is uniformly bounded on the set $0<r<s<t$ by an integrable
function not depending on $\varepsilon $. \ In fact, using the local
nondeterminism property \textbf{(LND)}, and Lemma \ \ref{lem2} in the
Appendix,\textbf{\ \ }we obtain the following lower bound for the
conditional variance $\sigma _{r,s,t}^{2}=\mathrm{Var}%
(B_{t}^{H,i}-B_{s}^{H,i}|\mathcal{F}_{r})$:%
\begin{equation}
\sigma _{r,s,t}^{2}\geq \mathrm{Var}(B_{t}^{H,i}-B_{s}^{H,i}|\mathcal{F}%
_{s})=\mathrm{Var}(B_{t}^{H,i}|\mathcal{F}_{s})\geq k_{2}(t-s)^{2H}.
\label{j1}
\end{equation}%
We can get rid off the factor $A_{r,t,s}^{i}$ in the expression (\ref{h5})
of $\Sigma _{\varepsilon }^{i}(r,t,s)$ using the inequality%
\begin{equation}
p_{t}(x)\leq C\frac{t^{-\frac{d}{2}+\frac{1}{2}}}{|x|}e^{-\frac{|x|^{2}}{4t}%
}\leq C\frac{t^{-\frac{d}{2}+\frac{1}{2}}}{|x|},  \label{j2}
\end{equation}%
for some constant $C>0$. In this way we obtain, using (\ref{j1}) and (\ref%
{j2})%
\begin{equation}
\left| \Sigma _{\varepsilon }^{i}(r,t,s)\right| \leq C\left( t-s\right)
^{-Hd-H}\left| K_{H}(t,r)-K_{H}(s,r)\right| ,  \label{j3a}
\end{equation}%
for some constant $C>0$, and by Lemma \ref{lem1} in the Appendix we obtain
that%
\begin{equation}
\int_{r}^{T}\int_{r}^{t}\left( t-s\right) ^{-Hd-H}\left|
K_{H}(t,r)-K_{H}(s,r)\right| dsdt\leq C(r^{\frac{1}{2}-H}\vee 1).  \label{j3}
\end{equation}%
By dominated convergence we deduce the convergence of the integrals%
\begin{equation*}
\lim_{\varepsilon \downarrow 0}\int_{r}^{T}\int_{r}^{t}\Sigma _{\varepsilon
}^{i}(r,t,s)dsdt=\int_{r}^{T}\int_{r}^{t}\Sigma ^{i}(r,t,s)dsdt
\end{equation*}%
for all $(r,\omega )\in \lbrack 0,T]\times \Omega $, and a second
application of the dominated convergence theorem yields that $%
\int_{r}^{T}\int_{r}^{t}\Sigma _{\varepsilon }^{i}(r,t,s)dsdt$ converges in
\ $L^{2}([0,T]\times \Omega )$ to $\int_{r}^{T}\int_{r}^{t}\Sigma
^{i}(r,t,s)dsdt$. This implies the convergence of $L_{\varepsilon }^{(1)}\ $%
\ to 
\begin{equation*}
-\sum_{i=1}^{d}\int_{0}^{T}\left( \int_{r}^{T}\int_{r}^{t}\Sigma
^{i}(r,t,s)dsdt\right) dW_{r}^{i}
\end{equation*}%
in $L^{2}(\Omega )$ as $\varepsilon $ tends to zero.

\textbf{Step 4} Consider now the case $s<r<t$. \ In this case the integral
of the term $\Sigma _{\varepsilon }^{i}(r,t,s)$ is not necessarily bounded,
and in order to show the convergence of \ $L_{\varepsilon }^{(2)}$ we will
prove uniform bounds in $\varepsilon $ for the expectation $E\left(
\int_{r}^{T}\int_{r}^{t}\ \left| \Sigma _{\varepsilon }^{i}(r,t,s)\right|
^{p}dsdt\right) $, for some $p>1$. We can write for $s<r<t$, using the first
inequality in (\ref{j2})%
\begin{eqnarray}
\left| \Sigma _{\varepsilon }^{i}(r,t,s)\right|  &\leq &\frac{\left|
A_{r,t,s}\right| }{\left( \varepsilon +\sigma _{r,s,t}^{2}\right) }%
p_{\varepsilon +\sigma _{r,s,t}^{2}}(A_{r,t,s})\left| K_{H}(t,r)\right|  
\notag \\
&=&(2\pi )^{-\frac{d}{2}}\frac{\left| A_{r,t,s}\right| }{\left( \varepsilon
+\sigma _{r,s,t}^{2}\right) ^{1+\frac{d}{2}}}\exp \left( -\frac{%
|A_{r,t,s}|^{2}}{2(\varepsilon +\sigma _{r,s,t}^{2})}\right) \left|
K_{H}(t,r)\right|   \notag \\
&\leq &C\left( \varepsilon +\sigma _{r,s,t}^{2}\right) ^{-\frac{d+1}{2}}\exp
\left( -\frac{|A_{r,t,s}|^{2}}{4(\varepsilon +\sigma _{r,s,t}^{2})}\right)
\left| K_{H}(t,r)\right| ,  \label{y3}
\end{eqnarray}%
for some constant $C>0$. If $s<r<t$, using the local nondeterminism property 
\textbf{(LND) \ }we obtain the following lower bound for the conditional
variance $\sigma _{r,s,t}^{2}:$%
\begin{equation}
\sigma _{r,s,t}^{2}=\mathrm{Var}(B_{t}^{H,i}-B_{s}^{H,i}|\mathcal{F}_{r})=%
\mathrm{Var}(B_{t}^{H,i}|\mathcal{F}_{r})\geq k_{2}(t-r)^{2H}.  \label{y2}
\end{equation}%
On the other hand, if $s<r<t$%
\begin{eqnarray}
\sigma _{r,s,t}^{2} &=&\mathrm{Var}(B_{t}^{H,i}-B_{s}^{H,i}|\mathcal{F}_{r})=%
\mathrm{Var}(B_{t}^{H,i}-B_{r}^{H,i}|\mathcal{F}_{r})  \notag \\
&\leq &\mathrm{Var}(B_{t}^{H,i}-B_{r}^{H,i})=(t-r)^{2H}.  \label{y1}
\end{eqnarray}%
Also we will make use of the estimate (see \cite{Hu2})%
\begin{equation}
\left| K_{H}(t,r)\right| \leq k_{3}(t-r)^{H-\frac{1}{2}}r^{\frac{1}{2}-H}.
\label{l3}
\end{equation}%
Substituting the estimates (\ref{y2}), (\ref{y1}) and (\ref{l3}) into \ (\ref%
{y3}) yields%
\begin{equation}
\left| \Sigma _{\varepsilon }^{i}(r,t,s)\right| \leq Cr^{\frac{1}{2}-H}\Psi
_{\varepsilon }(r,t,s),  \label{y5}
\end{equation}%
for some constant \ $C$, where%
\begin{equation}
\Psi _{\varepsilon }(r,t,s)=\left( \varepsilon +k_{2}\left( t-r\right)
^{2H}\right) ^{-\frac{d+1}{2}}(t-r)^{H-\frac{1}{2}}\exp \left( -\frac{%
|A_{r,t,s}|^{2}}{4(\varepsilon +(t-r)^{2H})}\right) .  \label{y4}
\end{equation}

Notice that if $Hd<\frac{1}{2}$, then $\left| \Sigma _{\varepsilon
}^{i}(r,t,s)\right| $ is uniformly bounded by the integrable function $Cr^{%
\frac{3}{2}-H}\ \left( t-r\right) ^{-Hd-\frac{1}{2}}$, and we can conclude
as in Step 3. For this reason, we can assume that $Hd\geq \frac{1}{2}$.

We claim that for some $p>1$, we have%
\begin{equation}
\sup_{\varepsilon >0}E\left( \int_{r}^{T}\int_{0}^{r}\Psi _{\varepsilon
}^{p}(r,t,s)dsdt\right) <\infty .  \label{eq1}
\end{equation}%
To show this estimate we first derive a lower bound for the expectation of $%
|A_{r,t,s}^{1}|^{2}=\left[ E(B_{t}^{H,1}-B_{s}^{H,1}|\mathcal{F}_{r})\right]
^{2}$. \ The main idea is to add and substract the term $B_{r}^{H,1}$, and
then neglect the expectation $E\left( \left( \left( E(B_{t}^{H,1}|\mathcal{F}%
_{r})-B_{r}^{H,1}\right) ^{2}\right) \right) $. This argument will be used
later to \ find a lower bound for the covariance matrix of the vector $%
\left( E(B_{t_{i}}^{H,1}-B_{s_{i}}^{H,1}|\mathcal{F}_{r}),1\leq i\leq
n\right) $. 
\begin{eqnarray*}
E\left( |A_{r,t,s}^{1}|^{2}\right)  &=&E\left( \left(
E(B_{t}^{H,1}-B_{s}^{H,1}|\mathcal{F}_{r})\right) ^{2}\right)  \\
&=&E\left( \left( E(B_{t}^{H,1}|\mathcal{F}_{r})-B_{r}^{H,1}\right)
^{2}\right)  \\
&&+2E\left( \left( E(B_{t}^{H,1}|\mathcal{F}_{r})-B_{r}^{H,1}\right) \left(
B_{r}^{H,1}-B_{s}^{H,1}\right) \right) +E\left( \left(
B_{r}^{H,1}-B_{s}^{H,1}\right) ^{2}\right)  \\
&\geq &2E\left( \left( B_{t}^{H,1}-B_{r}^{H,1}\right) \left(
B_{r}^{H,1}-B_{s}^{H,1}\right) \right) +E\left( \left(
B_{r}^{H,1}-B_{s}^{H,1}\right) ^{2}\right)  \\
&=&E\left( \left( B_{t}^{H,1}-B_{s}^{H,1}\right) ^{2}\right) -E\left( \left(
B_{t}^{H,1}-B_{r}^{H,1}\right) ^{2}\right)  \\
&=&(t-s)^{2H}-(t-r)^{2H}.
\end{eqnarray*}%
As a consequence, we obtain, assuming $p<2$ 
\begin{eqnarray*}
&&\!\!\!\!\!\!\!\!\!\!\!\!\!\!\!\!\!\!\!\!\!\!\!\!E\left( \exp \left( -\frac{%
p|A_{r,t,s}|^{2}}{4(\varepsilon +(t-r)^{2H})}\right) \right)  \\
&=&\left( 1+\frac{p}{2}(\varepsilon +(t-r)^{2H})^{-1}E\left(
|A_{r,t,s}^{1}|^{2}\right) \right) ^{-\frac{d}{2}} \\
&\leq &\left( 1+\frac{p}{2}(\varepsilon +(t-r)^{2H})^{-1}\left[
(t-s)^{2H}-(t-r)^{2H}\right] \right) ^{-\frac{d}{2}} \\
&=&(\varepsilon +(t-r)^{2H})^{\frac{d}{2}} \\
&&\times \left( \varepsilon +\left( 1-\frac{p}{2}\right) (t-r)^{2H}+\frac{p}{%
2}(t-s)^{2H}\right) ^{-\frac{d}{2}}.
\end{eqnarray*}%
Hence,%
\begin{eqnarray}
  && E\left( \exp \left( -\frac{%
p|A_{r,t,s}|^{2}}{4(\varepsilon +(t-r)^{2H})}\right) \right)   \notag  \\
 &&  \quad \leq
C(\varepsilon +(t-r)^{2H})^{\frac{d}{2}}(t-r)^{-2H\alpha }(t-s)^{-2H\beta },
\label{j8}
\end{eqnarray}%
where $\alpha +\beta =\frac{d}{2}$. Substituting (\ref{j8}) into (\ref{y4}) 
yields 
\begin{eqnarray*}
E\left( \int_{r}^{T}\int_{0}^{r}\Psi _{\varepsilon }^{p}(r,t,s)dsdt\right) 
&\leq &C\int_{r}^{T}\int_{0}^{r}\left( \varepsilon +\left( t-r\right)
^{2H}\right) ^{-\frac{d+1}{2}p+\frac{d}{2}-\alpha } \\
&&\times (t-r)^{\left( H-\frac{1}{2}\right) p}(t-s)^{-\beta 2H}dsdt \\
&\leq &C\int_{r}^{T}\int_{0}^{r}\left( t-r\right) ^{-pHd-\frac{p}{2}+2H\beta
}(t-s)^{-2H\beta }dsdt.
\end{eqnarray*}%
If $Hd>1$, we can choose $\beta $ such that $2H\beta >1$, and integrating in
the variable $s$, the above integral is bounded by%
\begin{equation*}
C\int_{r}^{T}\left( t-r\right) ^{-pHd-\frac{p}{2}+1}dt,
\end{equation*}%
which is finite it $p>1$ satisfyes $\left( Hd+\frac{1}{2}\right) p<2$ (this
is possible because $Hd+\frac{1}{2}<2$). \ If $Hd\leq 1$, we can choose $%
\beta $ such that $2H\beta =Hd-\delta $, for any $\delta >0$ , and we obtain
the bound%
\begin{equation*}
C\int_{r}^{T}\left( t-r\right) ^{-pHd-\frac{p}{2}+Hd-\delta }dt,
\end{equation*}%
which is again finite if $p>1$ is close to one, and $\delta >0$ is small
enough.

As a consequence, from (\ref{y5}) and (\ref{eq1}),  for any fixed $r\in \lbrack 0,T]$, the family of functions $\left\{ \ \Sigma _{\varepsilon
}^{i}(r,t,s),\varepsilon >0\right\} $, is uniformly integrable in $%
[r,T]\times \lbrack 0,r]$, so it converges in \ $L^{1}([r,T]\times \lbrack
0,r]) \times \Omega$ to \ $\Sigma ^{i}(r,t,s)$, for $i=1,\ldots ,d$. This implies the
convergence of the integrals%
\begin{equation*}
\lim_{\varepsilon \downarrow 0}\int_{r}^{T}\int_{0}^{r}\Sigma _{\varepsilon
}^{i}(r,t,s)dsdt=\int_{r}^{T}\int_{0}^{r}\Sigma ^{i}(r,t,s)dsdt,
\end{equation*}
for each fixed $r\in \lbrack 0,T]$ in $L^1(\Omega)$. \ 

Finally, we claim that this convergence also holds in \ $L^{2}([0,T]\times
\Omega )$, and this implies the convergence of $L_{\varepsilon }^{(2)}\ $\
to 
\begin{equation*}
-\sum_{i=1}^{d}\int_{0}^{T}\left( \int_{r}^{T}\int_{0}^{r}\Sigma
^{i}(r,t,s)dsdt\right) dW_{r}^{i}
\end{equation*}%
in $L^{2}(\Omega )$ as $\varepsilon $ tends to zero. To show the convergence
in $L^{2}([0,T]\times \Omega )$ of the integrals 
\begin{equation*}
Y_{\varepsilon }^{i}(r)=\int_{r}^{T}\int_{0}^{r}\Sigma _{\varepsilon
}^{i}(r,t,s)dsdt
\end{equation*}
it suffices to prove that%
\begin{equation}
\sup_{\varepsilon >0}\int_{0}^{T}E\left( \left| Y_{\varepsilon
}^{i}(r)\right| ^{p}\right) dr<\infty  \label{j7}
\end{equation}%
for all $i=1,\ldots ,d$ and for some $p>2$. The proof of (\ref{j7}) will be
the last step in the proof of this theorem.

\noindent \textbf{Step 5 }Suppose first that $Hd<1$. Then, from (\ref{y5})
we obtain%
\begin{equation*}
\int_{0}^{T}E\left( \left| Y_{\varepsilon }^{i}(r)\right| ^{p}\right) dr\leq
C\int_{0}^{T}E\left[ \left( \int_{r}^{T}\int_{0}^{r}\Psi _{\varepsilon
}(r,t,s)dsdt\right) ^{p}\right] r^{p\left( \frac{1}{2}-H\right) }dr.
\end{equation*}%
Using (\ref{y4}) and Minkowski's inequality yields%
\begin{eqnarray}
\left\| \int_{r}^{T}\int_{0}^{r}\Psi _{\varepsilon }(r,t,s)dsdt\right\| _{p}
&\leq &\int_{r}^{T}\int_{0}^{r}\left( \varepsilon +k_{2}\left( t-r\right)
^{2H}\right) ^{-\frac{d+1}{2}}(t-r)^{H-\frac{1}{2}}  \notag \\
&&\times \left\| \exp \left( -\frac{|A_{r,t,s}|^{2}}{4(\varepsilon
+(t-r)^{2H})}\right) \right\| _{p}dsdt,  \label{b7}
\end{eqnarray}%
and from (\ref{j8}), choosing $\beta =\frac{d}{2}$,  we get%
\begin{equation}
\left\| \exp \left( -\frac{|A_{r,t,s}|^{2}}{4(\varepsilon +(t-r)^{2H})}%
\right) \right\| _{p}\leq C(\varepsilon +(t-r)^{2H})^{\frac{d}{2p}}(t-s)^{-%
\frac{Hd}{p}}.  \label{b6}
\end{equation}%
Substituting (\ref{b6}) into (\ref{b7}) yields%
\begin{equation*}
\left\| \int_{r}^{T}\int_{0}^{r}\Psi _{\varepsilon }(r,t,s)dsdt\right\|
_{p}\leq C\int_{r}^{T}(t-r)^{-Hd-\frac{1}{2}+\frac{Hd}{p}}dr,
\end{equation*}%
which is finite if we choose $p>2$ such that $p<\frac{2Hd}{2Hd-1}$. Finally,
if $p\left( \frac{1}{2}-H\right) >-1$ we complete the proof of (\ref{j7}) in
the case $Hd<1$. 

In the case $Hd\geq 1$ we cannot apply the previous arguments, and  the
proof \ of (\ref{j7}) follows from the moment  estimates given in
Proposition \ref{prop2}.
\end{proof}

\bigskip

\noindent \textbf{Remark 1 }Theorem \ref{thm1} also provides an alternative
proof of the existence of the self-intersection local time in the case $H\in
\lbrack \frac{1}{d},\min  ( \frac{3}{2d},\frac{2}{d+1} ) )$, which
was proved by Hu and Nualart in \cite{HN} in the general case $Hd<\frac{3}{2}
$. Notice that for $d\geq 3$, the condition \ $H\in \lbrack \frac{1}{d},\min
 ( \frac{3}{2d},\frac{2}{d+1} ) )$ is equivalent to $1\leq Hd<\frac{%
3}{2}$, and for $d=2$ we require $H<\frac{2}{3}$, instead of the more
general condition $H<\frac{3}{4}$, that guarantees the existence of the
renormalized local time (see \cite{Ro} and \cite{HN}).

\medskip

The next Proposition contains the \ basic estimates on the moments of the
quadratic variation of the stochastic integral appearing in the
representation of the renormalized self-intersection local time.

\begin{proposition}
\label{prop2} Assume $1\leq Hd<\frac{3}{2}$. Set%
\begin{equation*}
\Lambda _{\varepsilon }(r)=\int_{r}^{T}\int_{0}^{r}\Psi _{\varepsilon
}(r,t,s)dsdt,
\end{equation*}%
where $\Psi _{\varepsilon }(r,t,s)$ has been defined in (\ref{y4}). Then,
for any integer $n\geq 1$,%
\begin{equation*}
E\left( \Lambda _{\varepsilon }^{n}(r)\right) \leq C^{n}(n!)^{\gamma },
\end{equation*}%
for some constant $C>0$, where%
\begin{equation*}
\gamma >\left( \frac{1}{2}+H\right) \left( d-\frac{1}{2H}\right) .
\end{equation*}
\end{proposition}

\begin{proof}
Set $g_{\varepsilon }(t-r)=\left( \varepsilon +k_{2}\left( t-r\right)
^{2H}\right) ^{-\frac{d+1}{2}}\left( t-r\right) ^{H-\frac{1}{2}}$. We have%
\begin{eqnarray}
E\left( \Lambda _{\varepsilon }^{n}(r)\right) &=&E\left[ \left(
\int_{r}^{T}\int_{0}^{r}g_{\varepsilon }(t-r)\exp \left( -\frac{%
|A_{r,s,t}|^{2}}{4(\varepsilon +(t-r)^{2H})}\right) dsdt\right) ^{n}\right] 
\notag \\
&=&n!\int_{[r,T]^{n}}\int_{S_{n}}\prod_{i=1}^{n}g_{\varepsilon }(t_{i}-r) 
\notag \\
&&\times \left( E\left( \exp \left( -\sum_{i=1}^{n}\frac{%
|A_{r,s_{i},t_{i}}^{1}|^{2}}{4(\varepsilon +(t_{i}-r)^{2H})}\right) \right)
\right) ^{d}dsdt\mathbf{,}  \label{a7}
\end{eqnarray}%
where $S_{n}=\{0<s_{1}<\cdots <s_{n}<r\}$, $s=(s_{1},\ldots ,s_{n})$ and $%
t=(t_{1},\ldots ,t_{n})$.

We denote by $Q$ the covariance matrix of the vector%
\begin{equation*}
\left( E(B_{t_{1}}^{H,1}-B_{s_{1}}^{H,1}|\mathcal{F}_{r}),\ldots
,E(B_{t_{n}}^{H,1}-B_{s_{n}}^{H,1}|\mathcal{F}_{r})\right) .
\end{equation*}%
Then, a well-known formula for Gaussian random variables implies that%
\begin{eqnarray}
E\left[ \exp \left( -\sum_{i=1}^{n}\frac{|A_{r,s_{i},t_{i}}^{1}|^{2}}{%
4(\varepsilon +(t_{i}-r)^{2H})}\right) \right]  &=&\det \left( I+\frac{1}{2}%
QD^{-1}\right) ^{-\frac{1}{2}}  \notag \\
&&\!\!\!\!\!\!\!\!\!\!\!\!\!\!\!\!\!\!\!\!\!\!\!\!=2^{\frac{n}{2}%
}\prod_{i=1}^{n}\sqrt{a_{i}}\det \left( 2D+Q\right) ^{-\frac{1}{2}},
\label{a8}
\end{eqnarray}%
where $D$ denotes the \ $n\times n$ diagonal matrix with entries $%
a_{i}=\varepsilon +$ $(t_{i}-r)^{2H}$. \ As in the computation of $E\left(
|A_{r,t,s}^{1}|^{2}\right) $, adding and substracting the term $B_{r}^{H,1}$
yields%
\begin{eqnarray*}
Q_{ij} &=&E\left( E(B_{t_{i}}^{H,1}-B_{s_{i}}^{H,1}|\mathcal{F}%
_{r})E(B_{t_{j}}^{H,1}-B_{s_{j}}^{H,1}|\mathcal{F}_{r})\right)  \\
&=&E\left( E(B_{t_{i}}^{H,1}-B_{r}^{H,1}|\mathcal{F}%
_{r})E(B_{t_{j}}^{H,1}-B_{r}^{H,1}|\mathcal{F}_{r})\right)  \\
&&+E\left( (B_{r}^{H,1}-B_{s_{i}}^{H,1})(B_{t_{j}}^{H,1}-B_{r}^{H,1})\right)
+E\left( (B_{t_{i}}^{H,1}-B_{r}^{H,1})(B_{r}^{H,1}-B_{s_{j}}^{H,1})\right) 
\\
&&+E\left( (B_{r}^{H,1}-B_{s_{i}})(B_{r}^{H,1}-B_{s_{j}}^{H,1})\right)  \\
&=&E\left( E(B_{t_{i}}^{H,1}-B_{r}^{H,1}|\mathcal{F}%
_{r})E(B_{t_{j}}^{H,1}-B_{r}^{H,1}|\mathcal{F}_{r})\right)  \\
&&-E\left( (B_{t_{i}}^{H,1}-B_{r}^{H,1})(B_{t_{j}}^{H,1}-B_{r}^{H,1})\right)
+E\left(
(B_{t_{i}}^{H,1}-B_{s_{i}}^{H,1})(B_{t_{j}}^{H,1}-B_{s_{j}}^{H,1})\right) .
\end{eqnarray*}%
Hence, we obtain%
\begin{equation*}
Q=R-N+M,
\end{equation*}%
where 
\begin{eqnarray*}
R_{ij} &=&E\left( E(B_{t_{i}}^{H,1}-B_{r}^{H,1}|\mathcal{F}%
_{r})E(B_{t_{j}}^{H,1}-B_{r}^{H,1}|\mathcal{F}_{r})\right) , \\
M_{ij} &=&E\left(
(B_{t_{i}}^{H,1}-B_{s_{i}}^{H,1})(B_{t_{j}}^{H,1}-B_{s_{j}}^{H,1})\right) ,
\\
N_{ij} &=&E\left(
(B_{t_{i}}^{H,1}-B_{r}^{H,1})(B_{t_{j}}^{H,1}-B_{r}^{H,1})\right) .
\end{eqnarray*}%
\ All these matrices are nonnegative definite. The main idea will be to get
rid off the matrix $R$, and control the matrix $N$ $\ $by its diagonal
elements which are%
\begin{equation*}
N_{ii}=(t_{i}-r)^{2H}.
\end{equation*}%
Indeed, the matrix $N$ is nonnegative definite and, hence, it safisties the
inequality%
\begin{equation}
N\leq nD_{N},  \label{u1}
\end{equation}
where $D_{N}$ is a diagonal matrix whose entries are $N_{ii}$. Therefore,%
\begin{equation*}
Q\geq -N+M\geq -nD_{N}+M,
\end{equation*}%
and for any $1\leq \delta <2$, we can write 
\begin{equation}
\det (2D+Q)\geq \det (2D+\frac{2-\delta }{n}Q)\leq \det (2D-(2-\delta )D_{N}+%
\frac{2-\delta }{n}M).  \label{a9}
\end{equation}%
The entries of the diagonal matrix $D_{1}=2D-(2-\delta )D_{N}$ are the
positive numbers 
\begin{equation*}
2\varepsilon +\delta (t_{i}-r)^{2H}>0.
\end{equation*}%
From (\ref{a7}), (\ref{a8}) and (\ref{a9}) we obtain%
\begin{eqnarray*}
E\left( \Lambda _{\varepsilon }^{n}(r)\right)  &\leq &2^{\frac{nd}{2}%
}n!\int_{[r,T]^{n}}\int_{S_{n}}\prod_{i=1}^{n}\left( g_{\varepsilon
}(t_{i}-r)a_{i}^{\frac{d}{2}}\right)  \\
&&\times \det (D_{1}+\frac{2-\delta }{n}M)^{-\frac{d}{2}}dsdt\mathbf{.}
\end{eqnarray*}%
We have%
\begin{equation*}
\det (D_{1}+\frac{2-\delta }{n}M)^{-\frac{d}{2}}\leq \left( \frac{n}{%
2-\delta }\right) ^{n\beta }\left( \det D_{1}\right) ^{-\alpha }\left( \det
M\right) ^{-\beta },
\end{equation*}%
where $\alpha +\beta =\frac{d}{2}$. Hence,%
\begin{eqnarray*}
E\left( \Lambda _{\varepsilon }^{n}(r)\right)  &\leq &\ \left( \frac{n}{%
2-\delta }\right) ^{n\beta }2^{\frac{nd}{2}}n!\int_{[r,T]^{n}}\int_{S_{n}}%
\prod_{i=1}^{n}\left( g_{\varepsilon }(t_{i}-r)a_{i}^{\frac{d}{2}}\left(
2\varepsilon +\delta (t_{i}-r)^{2H}\right) ^{-\alpha }\right)  \\
&&\times (\det M)^{-\beta }dsdt\mathbf{.}
\end{eqnarray*}%
Then,%
\begin{eqnarray*}
&&g_{\varepsilon }(t_{i}-r)a_{i}^{\frac{d}{2}}\left( 2\varepsilon
+2(t_{i}-r)^{2H}\right) ^{-\alpha } \\
&=&\left( \varepsilon +k_{2}\left( t_{i}-r\right) ^{2H}\right) ^{-\frac{d+1}{%
2}}\left( t_{i}-r\right) ^{H-\frac{1}{2}}\left( \varepsilon
+(t_{i}-r)^{2H}\right) ^{\frac{d}{2}}\left( 2\varepsilon
+2(t_{i}-r)^{2H}\right) ^{-\alpha } \\
&\leq &C(t_{i}-r)^{-\frac{1}{2}-2H\alpha },
\end{eqnarray*}%
for some constant $C>0$. Thus%
\begin{equation}
E\left( \Lambda _{\varepsilon }^{n}(r)\right) \leq C^{n}n^{\beta
n}n!\int_{[r,T]^{n}}\int_{S_{n}}\prod_{i=1}^{n}(t_{i}-r)^{-\frac{1}{2}%
-2H\alpha }(\det M)^{-\beta }dsdt\mathbf{,}  \label{a2}
\end{equation}%
for some constant $C>0$.

Applying \ Lemma \ref{lem2} in the Appendix and the local nondeterminism
property of the fractional Brownian motion we obtain%
\begin{align}
\det M& =\mathrm{Var}(B_{t_{n}}-B_{s_{n}})\mathrm{Var}(B_{t_{{n-1}}}-B_{s_{{%
n-1}}}|B_{t_{{n}}}-B_{s_{n}})  \notag \\
& \times \cdots \times \mathrm{Var}(B_{t_{1}}-B_{s_{1}}|B_{t_{2}}-B_{s_{2}},%
\dots ,B_{t_{n}}-B_{_{n}})  \notag \\
& =(t_{n}-s_{n})^{2H}\mathrm{Var}(B_{s_{{n-1}}}|B_{t_{{n-1}}},B_{t_{{n}%
}},B_{s_{{n}}})  \notag \\
& \times \cdots \times \mathrm{Var}(B_{s_{1}}|B_{t_{1}},\dots
,B_{t_{n}},B_{s_{1}},\ldots ,B_{s_{n-1}})  \notag \\
\geq & k_{2}^{n-1}(r-s_{n})^{2H}\left( (s_{n}-s_{n-1})\wedge s_{n-1}\right)
^{2H}\cdots \left( (s_{2}-s_{1})\wedge s_{1}\right) ^{2H}.  \label{a1}
\end{align}%
Substituting (\ref{a1}) into (\ref{a2}), and choosing $\alpha $ such that $%
\alpha <\frac{1}{4H}$ (this is possible because $Hd\geq 1$) yields%
\begin{equation*}
E\left( \Lambda _{\varepsilon }^{n}(r)\right) \leq C^{n}n^{\beta
n}n!\int_{S_{n}}\left[ (r-s_{n})\left( (s_{n}-s_{n-1})\wedge s_{n-1}\right)
\cdots \left( (s_{2}-s_{1})\wedge s_{1}\right) \right] ^{-2\beta H}ds\mathbf{%
.}
\end{equation*}%
Finally, by Lemma \ref{lem3} in the Appendix we obtain 
\begin{equation*}
E\left( \Lambda _{\varepsilon }^{n}(r)\right) \leq \ \frac{C^{n}n^{\beta n}n!%
}{\Gamma (n(1-2H\beta )+1)}.
\end{equation*}%
Notice that $\beta =\frac{d}{2}-\alpha >\frac{d}{2}-\frac{1}{4H}$. And
hence, 
\begin{equation*}
E\left( \Lambda _{\varepsilon }^{n}(r)\right) \leq C^{n}(n!^{\beta +2H\beta
},
\end{equation*}%
where%
\begin{equation*}
\beta (1+2H)>\frac{d}{2}-\frac{1}{4H}+Hd-\frac{1}{2}=\left( \frac{1}{2}%
+H\right) \left( d-\frac{1}{2H}\right) .
\end{equation*}

This concludes the proof.
\end{proof}

Using the above proposition we can deduce the following integrability
results for the renormalized self-intersection local time.

\begin{theorem}
Assume $\frac{1}{d}\leq H<\min \left( \frac{3}{2d},\frac{2}{d+1}\right) $.
For any integer $p<\frac{1}{2}\left[ \left( \frac{1}{2}+H\right) \left( d-%
\frac{1}{2H}\right) \right] ^{-1}$ we have%
\begin{equation*}
E(\exp |\widetilde{L}|^{p})<\infty .
\end{equation*}
\end{theorem}

\begin{proof}
Taking into account Lemma \ref{relation} in the Appendix, it suffices to
show that%
\begin{equation*}
E\left( \exp \left\langle \widetilde{L}\right\rangle ^{p}\right) <\infty ,
\end{equation*}%
where%
\begin{equation*}
\left\langle \widetilde{L}\right\rangle =\sum_{i=1}^{d}\int_{0}^{T}\left(
\int_{r}^{T}\int_{0}^{t}\Sigma ^{i}(r,t,s)dsdt\right) ^{2}dr.
\end{equation*}%
As in the proof of \ Theorem \ref{thm1} we make the decomposition%
\begin{equation*}
\int_{r}^{T}\int_{0}^{t}\Sigma
^{i}(r,t,s)dsdt=\int_{r}^{T}\int_{r}^{t}\Sigma
^{i}(r,t,s)dsdt+\int_{r}^{T}\int_{0}^{r}\Sigma ^{i}(r,t,s)dsdt.
\end{equation*}%
From (\ref{j3a}) and (\ref{j3}) we know that%
\begin{equation*}
\left| \int_{r}^{T}\int_{r}^{t}\Sigma ^{i}(r,t,s)dsdt\right| \leq C(r^{\frac{%
1}{2}-H}\vee 1).
\end{equation*}%
Therefore, applying Fatou's lemma and the estimate (\ref{y5}) yields 
\begin{eqnarray*}
E(\exp \left\langle \widetilde{L}\right\rangle ^{p}) &\leq &CE\left( \exp
\left( \left| \sum_{i=1}^{d}\int_{0}^{T}\left(
\int_{r}^{T}\int_{0}^{r}\Sigma ^{i}(r,t,s)dsdt\right) ^{2}dr\right|
^{p}\right) \right)  \\
&\leq &C\lim \inf_{\varepsilon \downarrow 0}E\left( \exp \left( \left|
\sum_{i=1}^{d}\int_{0}^{T}\left( \int_{r}^{T}\int_{0}^{r}\Sigma
_{\varepsilon }^{i}(r,t,s)dsdt\right) ^{2}dr\right| ^{p}\right) \right)  \\
&\leq &C\lim \inf_{\varepsilon \downarrow 0}E\left( \exp \left( C\left|
\int_{0}^{T}r^{1-2H}\left( \int_{r}^{T}\int_{0}^{r}\Psi _{\varepsilon
}(r,t,s)dsdt\right) ^{2}dr\right| ^{p}\right) \right) .
\end{eqnarray*}%
Applying H\"{o}lder and Jensen inequalities we obtain 
\begin{eqnarray*}
E(\exp \left\langle \widetilde{L}\right\rangle ^{p}) &\leq &C\lim
\inf_{\varepsilon \downarrow 0}E\left( \exp \left(
C\int_{0}^{T}r^{1-2H}\left( \int_{r}^{T}\int_{0}^{r}\Psi _{\varepsilon
}(r,t,s)dsdt\right) ^{2p}dr\right) \right)  \\
&\leq &C\lim \inf_{\varepsilon \downarrow 0}\int_{0}^{T}r^{1-2H}E\left( \exp
\left( C\left( \int_{r}^{T}\int_{0}^{r}\Psi _{\varepsilon
}(r,t,s)dsdt\right) ^{2p}\right) \right) dr.
\end{eqnarray*}%
Finally,%
\begin{eqnarray*}
&&E\left( \exp \left( C\left( \int_{r}^{T}\int_{0}^{r}\Psi _{\varepsilon
}(r,t,s)dsdt\right) ^{2p}\right) \right)  \\
&=&\sum_{n=1}^{\infty }\frac{C^{n}}{n!}E\left( \left(
\int_{r}^{T}\int_{0}^{r}\Psi _{\varepsilon }(r,t,s)dsdt\right) ^{2np}\right) 
\\
&\leq &\sum_{n=1}^{\infty }\frac{C^{n}}{n!}(\left( \left[ 2np\right]
+1\right) !)^{\gamma },
\end{eqnarray*}%
and it suffices to apply Proposition \ref{prop2} to conclude the proof.
\end{proof}

\medskip 

\noindent
\noindent \textbf{Remark 2 }The exponent $p_{0}=\frac{1}{2}\left[ \left( 
\frac{1}{2}+H\right) \left( d-\frac{1}{2H}\right) \right] ^{-1}$ is not
optimal. For instance, if $Hd=1$, then  $p_{0}=\frac{2H}{1+2H}$ and we know
that for $Hd<1$, then $p_{0}=\frac{1}{Hd}$. In particular,  if $H=\frac{1}{2}
$ and $d=2$ we obtain $p_{0}=\frac{1}{2}$, and we know that in this case the
critical exponent is $p_{0}=1$. \  The lack of optimality is due to the
factor $n$ in the estimation of the \ positive definite matrix $N$ $\ $by
its diagonal elements given in (\ref{u1}). Without this factor $n$ we would
get the critical exponent $\frac{1}{2Hd-1}$, but our method does not allow
to get this value.

\medskip \noindent \textbf{Remark 3} In the case of the planar Brownian
motion $B=\{B_{t},t\geq 0\}$ (that is, $d=2$, and $H=\frac{1}{2}$), formula
\ (\ref{h1}) yields%
\begin{equation}
\widetilde{L}=-\frac{1}{2\pi }\sum_{i=1}^{2}\int_{0}^{T}\left(
\int_{r}^{T}\int_{0}^{r}\frac{B_{r}^{i}-B_{s}^{i}}{(t-r)^{2}}\exp \left( -%
\frac{|B_{r}-B_{s}|^{2}}{2(t-r)}\right) dsdt\right) dB_{r}^{i}.  \label{t6}
\end{equation}%
The quadratic variation of this stochastic integral is%
\begin{eqnarray*}
\left\langle \widetilde{L}\right\rangle  &=&\frac{1}{4\pi ^{2}}%
\sum_{i=1}^{2}\int_{0}^{T}\left( \int_{r}^{T}\int_{0}^{r}\frac{%
B_{r}^{i}-B_{s}^{i}}{(t-r)^{2}}\exp \left( -\frac{|B_{r}-B_{s}|^{2}}{2(t-r)}%
\right) dsdt\right) ^{2}dr \\
&\leq &\frac{1}{4\pi ^{2}}\int_{0}^{T}\left( \int_{r}^{T}\int_{0}^{r}\frac{%
\left| B_{r}-B_{s}\right| }{(t-r)^{2}}\exp \left( -\frac{|B_{r}-B_{s}|^{2}}{%
2(t-r)}\right) dsdt\right) ^{2}dr \\
&=&\frac{1}{\pi ^{2}}\int_{0}^{T}\left( \int_{0}^{r}\frac{1}{\left|
B_{r}-B_{s}\right| }\exp \left( -\frac{|B_{r}-B_{s}|^{2}}{2(T-r)}\right)
ds\right) ^{2}dr \\
&\leq &\frac{1}{\pi ^{2}}\int_{0}^{T}\left( \int_{0}^{r}\frac{ds}{\left|
B_{r}-B_{s}\right| }\right) ^{2}dr.
\end{eqnarray*}%
From It\^{o}'s calculus we know that%
\begin{equation*}
\int_{0}^{r}\frac{ds}{\left| B_{r}-B_{s}\right| }=\frac{1}{d-1}\left(
X_{r}-b_{r}\right) ,
\end{equation*}%
where $X_{r}$ has the law of the modulus of a $d$-dimensional Brownian
motion at time $r$ (Bessel process), and $b_{r}$ has a normal $N(0,r)$ law.
We can write%
\begin{equation*}
\exp \left( \lambda \left\langle \widetilde{L}\right\rangle \right) \leq 
\frac{1}{T}\int_{0}^{T}\exp \left( \frac{T\lambda }{\pi ^{2}}\left(
\int_{0}^{r}\frac{ds}{\left| B_{r}-B_{s}\right| }\right) ^{2}\right) dr,
\end{equation*}%
which clearly imply the existence of some $\lambda _{0}$ such that $E$ $%
\left( \exp \left( \lambda \left\langle \widetilde{L}\right\rangle \right)
\right) <\infty $ for all $\lambda <\lambda _{0}$. From Lemma \ \ref%
{relation} we get that \ there exists \ $\beta _{0}$ such that $E$ $\left(
\exp \left( \beta \left| \widetilde{L}\right| \right) \right) <\infty $ for
all $\beta <\beta _{0}$. This method does not allows us to obtain the
critical exponent, just the existence of exponential moments.

\medskip \noindent \textbf{Remark 4 }$\ $The above results remain true if we
replace the \ fractional Brownian motion with Hurst paramter $H$, by an
arbitrary \ centered Gaussian process of the form (\ref{g1}) satisfying the
local nondeterminism property \textbf{(LND)} and following properties:

\begin{itemize}
\item[(C1)] For any $s,t\in \lbrack 0,T]$, \ $s<t$, there exist \ constants $%
k_{3}$ and $k_{4}$ such that 
\begin{equation*}
\ k_{3}(t-s)^{2H}\leq E(|B_{t}^{i}-B_{s}^{i}|^{2})\leq k_{4}(t-s)^{2H}.
\end{equation*}

\item[(C2)] The kernel $K(t,s)$ satisfies the estimates%
\begin{equation*}
\left| K(t,s)\right| \leq k_{5}(t-s)^{H-\frac{1}{2}}s^{\frac{1}{2}-H}\text{,}
\end{equation*}
\ for all $s<t$, and%
\begin{equation*}
\int_{r}^{T}\int_{r}^{t}\left( t-s\right) ^{-Hd-H}\left|
K(t,r)-K(s,r)\right| dsdt\leq \ \psi (r),
\end{equation*}%
where $\int_{0}^{T}$ $\psi (r)^{2}dr<\infty $.
\end{itemize}

\setcounter{equation}{0}

\section{Appendix}

In this Appendix we will first state and prove some elementary lemmas. The
first one is well-known.

\begin{lemma}
\label{lem2} Suppose that $\mathcal{G}_{1}$ $\subset \mathcal{G}_{2}$ are
two $\sigma $-fields contained in $\mathcal{F}$. Then, for any square
integrable random variable $F$ we have%
\begin{equation*}
\mathrm{Var}(F|\mathcal{G}_{1})\geq \mathrm{Var}(F|\mathcal{G}_{2}).
\end{equation*}
\end{lemma}

Let \ $M=\{M_{t},t\geq 0\}$ be a continuous local martingale such that $%
M_{0}=0$. Then, the following maximal exponential inequality is well-known%
\begin{equation*}
P\left( \sup_{\substack{ 0\leq t\leq T}}|M_{t}|\geq \delta ,\langle M\rangle
_{T}<\rho \right) \leq 2\exp \left( -\frac{\delta ^{2}}{2\rho }\right) .
\end{equation*}%
As a consequence of this inequality we can obtain exponential moments for $%
M_{T}$ from exponential moments of the quadratic variation \ $\langle
M\rangle _{T}$

\begin{lemma}
\label{relation} Suppose that for some $\alpha >0$ and $p\in (0,1]$ we have $%
E(e^{\alpha \langle M\rangle _{T}^{p}})<\infty $. Then,

\begin{itemize}
\item[(i)]  if $p=1$, for any $\lambda <\sqrt{\dfrac{\alpha }{2}}$, $%
E(e^{\lambda |M_{T}|})<\infty $, and

\item[(ii)]  if $p<1$, $E(e^{\lambda |M_{T}|^{p}})<\infty $ for all $\lambda
>0$.
\end{itemize}
\end{lemma}

\begin{proof}
Set $X=|M_{T}|^{p}$. For any constant \ $c>0$ we can write 
\begin{eqnarray}
E(e^{\lambda X}) &=&\int_{0}^{\infty }P(X\geq y)\lambda e^{\lambda y}dy 
\notag \\
&=&\int_{0}^{\infty }\left[ P(X\geq y,\langle M\rangle _{T}^{p}<cy)+P(X\geq
y,\langle M\rangle _{T}^{p}\geq cy)\right] \lambda e^{\lambda y}dy  \notag \\
&\leq &\int_{0}^{\infty }2\exp \left( -\frac{y^{\frac{1}{p}}}{2c^{\frac{1}{p}%
}}\right) \lambda e^{\lambda y}dy+\int_{0}^{\infty }P\left( \dfrac{\langle
M\rangle _{T}^{p}}{c}\geq y\right) \lambda e^{\lambda y}dy  \notag \\
&=&\int_{0}^{\infty }2\lambda \exp \left( \lambda y-\frac{y^{\frac{1}{p}}}{%
2c^{\frac{1}{p}}}\right) dy+E(e^{\tfrac{\lambda }{c}\langle M\rangle
_{T}^{p}}).  \notag
\end{eqnarray}%
Then it suffices to choose $c=\frac{\lambda }{\alpha }$ to complete the
proof.
\end{proof}

\bigskip The next two results are technical lemmas used in the paper.

\begin{lemma}
\label{lem1} Suppose that $H<\min (\frac{2}{d+1},\frac{3}{2d})$. Then, we
have%
\begin{equation*}
\int_{r}^{T}\int_{r}^{t}\left( t-s\right) ^{-Hd-H}\left|
K_{H}(t,r)-K_{H}(s,r)\right| dsdt\leq C\left( r^{\frac{1}{2}-H}\vee 1\right)
,
\end{equation*}%
for some constant $C$.
\end{lemma}

\begin{proof}
We know that%
\begin{equation*}
\frac{\partial K_{H}}{\partial t}(t,s)=c_{H}\left( H-\frac{1}{2}\right)
\left( \frac{t}{s}\right) ^{H-\frac{1}{2}}(t-s)^{H-\frac{3}{2}}.
\end{equation*}%
Then%
\begin{eqnarray*}
&&I:=\int_{r}^{T}\int_{r}^{t}\left( t-s\right) ^{-Hd-H}\left|
K_{H}(t,r)-K_{H}(s,r)\right| dsdt \\
&\leq &C\int_{r}^{T}\int_{r}^{t}\int_{s}^{t}\left( t-s\right) ^{-Hd-H}\left( 
\frac{\theta }{r}\right) ^{H-\frac{1}{2}}(\theta -r)^{H-\frac{3}{2}}d\theta
dsdt.
\end{eqnarray*}%
If $H<\frac{1}{2}$, then, $\left( \frac{\theta }{r}\right) ^{H-\frac{1}{2}%
}\leq 1$, and if $H>\frac{1}{2}$, then $\left( \frac{\theta }{r}\right) ^{H-%
\frac{1}{2}}\leq Cr^{\frac{1}{2}-H}$. Hence, the above integral is bounded by%
\begin{equation*}
C(r^{\frac{1}{2}-H}\vee 1)\int_{r}^{T}\int_{r}^{t}\int_{s}^{t}\left(
t-s\right) ^{-Hd-H}(\theta -r)^{H-\frac{3}{2}}d\theta dsdt.
\end{equation*}%
From the decomposition%
\begin{eqnarray*}
\frac{3}{2}-H &=&\alpha +\beta , \\
Hd+H &=&\gamma +\delta ,
\end{eqnarray*}%
we obtain%
\begin{eqnarray*}
&&\int_{r}^{T}\int_{r}^{t}\int_{s}^{t}\left( t-s\right) ^{-Hd-H}(\theta
-r)^{H-\frac{3}{2}}d\theta dsdt \\
&=&\int_{r}^{T}\int_{r}^{t}\int_{s}^{t}\left( s-r\right) ^{-\alpha }(\theta
-s)^{-\beta -\gamma }(t-\theta )^{-\delta }d\theta dsdt.
\end{eqnarray*}%
Finally, it suffices to show the parameters $\alpha $, $\beta $, $\gamma $
and $\delta $ in such a way that $\alpha <1$, $\delta <1$ and $\beta +\gamma
<1$. This leads to the condition%
\begin{equation*}
\frac{1}{2}+Hd<\min (1,\frac{3}{2}-H)+\min (1,Hd+H),
\end{equation*}%
which is satisfied if $H<\min (\frac{2}{d+1},\frac{3}{2d})$.  
\end{proof}

\begin{lemma}
\label{lem3} Let $a<1$. Fix an interval $[0,T]$. For each integer $n\geq 1$
we have%
\begin{eqnarray}
&& \int_{\Delta _{n}(T)}\left[ \left( (T-s_{n})\wedge s_{n}\right) \left(
(s_{n}-s_{n-1})\wedge s_{n-1}\right) \cdots \left( (s_{2}-s_{1})\wedge
s_{1}\right) \right] ^{-a}ds  \notag \\
&& \quad \le \frac{T^{n(1-a)}}{\Gamma (n(1-a)+1)}C^{n}\mathbf{,}  \label{a3}
\end{eqnarray}%
where $\Delta _{n}(T)=\{0<s_{1}<\cdots <s_{n}<T\}$
\end{lemma}

\begin{proof}
We proceed by induction on $n$. \ For $n=1$ we can write%
\begin{eqnarray*}
\int_{0}^{T}\left( (T-s_{1})\wedge s_{1}\right) ^{-a}ds_{1} &=&\int_{0}^{%
\frac{T}{2}}s_{1}^{-a}ds_{1}+\int_{\frac{T}{2}}^{T}(T-s_{1})^{-a}ds_{1} \\
&=&\frac{2}{1-a}\left( \frac{T}{2}\right) ^{1-a},
\end{eqnarray*}%
which implies (\ref{a3}) with $C=\frac{\Gamma (2-a)}{1-a}2^{a}$.

Suppose that the result holds for $n-1$. Then, 
\begin{eqnarray*}
I_{n} &=&\int_{\Delta _{n}(T)}\left[ \left( (T-s_{n})\wedge s_{n}\right)
\left( (s_{n}-s_{n-1})\wedge s_{n-1}\right) \cdots \left(
(s_{2}-s_{1})\wedge s_{1}\right) \right] ^{-a}ds \\
&=&\int_{0}^{T}\left( (T-s_{n})\wedge s_{n}\right) ^{-a} \\
&&\times \left( \int_{\Delta _{n-1}(s_{n})}\left[ \left(
(s_{n}-s_{n-1})\wedge s_{n-1}\right) \cdots \left( (s_{2}-s_{1})\wedge
s_{1}\right) \right] ^{-a}ds_{1}\cdots ds_{n-1}\right) ds_{n}.
\end{eqnarray*}%
By the induction hypothesis we can write%
\begin{eqnarray*}
I_{n} &\leq &\frac{C^{n-1}}{\Gamma (n-a)}\int_{0}^{T}\left( (T-s_{n})\wedge
s_{n}\right) ^{-a}s_{n}^{(n-1)(1-a)}ds_{n} \\
&=&\frac{C^{n-1}}{\Gamma ((n-1)(1-a)+1)} \\
&&\times \left( \int_{0}^{\frac{T}{2}}s_{n}^{(n-1)(1-a)-a}ds_{n}+\int_{\frac{%
T}{2}}^{T}(T-s_{n})^{-a}s_{n}^{(n-1)(1-a)}ds_{n}\right)  \\
&\leq &\frac{C^{n-1}}{\Gamma (n(1-a)+a)} \\
&&\times \left( \frac{1}{n(1-a)}\left( \frac{T}{2}\right)
^{n(1-a)}+T^{n(1-a)}\int_{0}^{1}(1-x)^{-a}x^{(n-1)(1-a)}dx\right)  \\
&\leq &\frac{T^{n(1-a)}C^{n-1}}{\Gamma (n(1-a)+a)}\left( \frac{1}{n(1-a)}+%
\frac{\Gamma (1-a)\Gamma ((n-1)(1-a)+1)}{\Gamma (n(1-a)+1)}\right)  \\
&=&T^{n(1-a)}C^{n-1}\left( \frac{1}{n(1-a)\Gamma (n(1-a)+a)}+\frac{\Gamma
(1-a)}{\Gamma (n(1-a)+1)}\right) .
\end{eqnarray*}%
Using the relation $\Gamma (n+1)=n\Gamma (n)$ we obtain%
\begin{equation*}
n(1-a)\Gamma (n(1-a)+a)\geq n(1-a)\Gamma (n(1-a))=\Gamma (n(1-a)+1),
\end{equation*}%
and, as a consequence%
\begin{equation*}
I_{n}\leq T^{n(1-a)}C^{n-1}\left( 1+\Gamma (1-a)\right) \frac{1}{\Gamma
(n(1-a)+1)},
\end{equation*}%
and it suffices to take $C\geq \max \left( \frac{\Gamma (2-a)}{1-a}%
2^{a},1+\Gamma (1-a)\right) .$
\end{proof}

\end{document}